\documentclass[12pt]{amsart}%
\usepackage{graphicx}
\usepackage{amscd}
\usepackage{amsmath}
\usepackage{amsfonts}
\usepackage{amssymb}%
\providecommand{\U}[1]{\protect\rule{.1in}{.1in}}
\theoremstyle{plain}

\newtheorem{corollary}{Corollary}[section]

\newtheorem{definition}{Definition}[section]

\newtheorem{lemma}{Lemma}[section]

\newtheorem{theorem}{Theorem}
\theoremstyle{definition}

\newtheorem{remark}{Remark}[section]

\numberwithin{equation}{section}
\numberwithin{theorem}{section}
\begin{document}
\title[Bessel-type functions]{Fourth-order Bessel-type special functions:\\a survey}
\author{W.N. Everitt}
\address{W.N. Everitt, School of Mathematics and Statistics, University of Birmingham,
Edgbaston, Birmingham B15 2TT, England, UK}
\email{w.n.everitt@bham.ac.uk}
\date{25 August 2005 (File: C:%
$\backslash$%
Swp50%
$\backslash$%
Bessel%
$\backslash$%
munich9.tex)}
\dedicatory{This paper is dedicated to the memory and achievements of\\George Neville Watson (1886 to 1965)}
\begin{abstract}
This survey paper reports on the properties of the fourth-order Bessel-type
linear ordinary differential equation, on the generated self-adjoint
differential operators in two associated Hilbert function spaces, and on the
generalisation of the classical Hankel integral transform.

These results are based upon the properties of the classical Bessel and
Laguerre second-order differential equations, and on the fourth-order
Laguerre-type differential equation. From these differential equations and
their solutions, limit processes yield the fourth-order Bessel-type functions
and the associated differential equation.

\end{abstract}
\maketitle
\tableofcontents

\section{Introduction\label{sec1}}

This survey paper is based on joint work with the following named colleagues:

\begin{center}
Jyoti Das, University of Calcutta, India

D.B. Hinton, University of Tennessee, USA

H. Kalf, University of Munich, Germany

L.L. Littlejohn, Utah State University, USA

C. Markett, Technical University of Aachen, Germany

M. Plum, University of Karlsruhe, Germany

M. van Hoeij, Florida State University, USA
\end{center}

\section{History\label{sec2}}

We see below that the structured definition of the general-even order
Bessel-type special functions is dependent upon the Jacobi and Laguerre
classical orthogonal polynomials, and the Jacobi-type and Laguerre-type
orthogonal polynomials.

These latter orthogonal polynomials were first defined by H.L. Krall in 1940,
see \cite{HK} and \cite{HK1}, and later studied in detail by A.M. Krall in
1981, see \cite{AK}, and by Koornwinder in 1984, see \cite{THK}. In this
respect see the two survey papers, \cite{EL1} of 1990 and \cite{PAT} of 1999.

The Bessel-type special functions of general even-order were introduced by
Everitt and Markett in 1994, see \cite{EM}.

The properties of the fourth-order Bessel-type functions have been studied by
the present author and the seven colleagues named in Section \ref{sec1} above,
in the papers \cite{DEHLM}, \cite{EKLM} and \cite{EKLM1}.

\section{The fourth-order differential equation\label{sec3}}

The fourth-order Bessel-type differential equation takes the form%
\begin{equation}
(xy^{\prime\prime}(x))^{\prime\prime}-((9x^{-1}+8M^{-1}x)y^{\prime
}(x))^{\prime}=\Lambda xy(x)\ \text{for all}\ x\in(0,\infty) \label{eq3.1}%
\end{equation}
where $M\in(0,\infty)$ is a positive parameter and $\Lambda\in\mathbb{C},$ the
complex field, is a spectral parameter. The differential equation
(\ref{eq3.1}) is derived in the paper \cite[Section 1, (1.10a)]{EM}, by
Everitt and Markett.

This linear, ordinary differential equation on the interval $(0,\infty
)\subset\mathbb{R},$ the real field, is written in Lagrange symmetric
(formally self-adjoint) form, or equivalently Naimark form, see \cite[Chapter
V]{MAN}.

The structured Bessel-type functions of all even orders, and their associated
linear differential equations, were introduced in the paper \cite[Section
1]{EM} through linear combinations of, and limit processes applied to, the
Laguerre and Laguerre-type orthogonal polynomials, and to the classical Bessel
functions. This process is best illustrated through the following diagram, see
\cite[Section 1, Page 328]{EM} (for the first two lines of this table see the
earlier work of Koornwinder \cite{THK} and Markett \cite{CM}):%
\begin{equation}
\left\{
\begin{array}
[c]{ccc}%
\begin{array}
[c]{c}%
\text{Jacobi polynomials}\\
k(\alpha,\beta)(1-x)^{\alpha}(1+x)^{\beta}%
\end{array}
& \longrightarrow &
\begin{array}
[c]{c}%
\text{Jacobi-type polynomials}\\
k(\alpha,\beta)(1-x)^{\alpha}(1+x)^{\beta}\\
+M\delta(x+1)+N\delta(x-1)
\end{array}
\\
\downarrow &  & \downarrow\\%
\begin{array}
[c]{c}%
\text{Laguerre polynomials}\\
k(\alpha)x^{\alpha}\exp(-x)
\end{array}
& \longrightarrow &
\begin{array}
[c]{c}%
\text{Laguerre-type polynomials}\\
k(\alpha)x^{\alpha}\exp(-x)+N\delta(x)
\end{array}
\\
\downarrow &  & \downarrow\\%
\begin{array}
[c]{c}%
\text{Bessel functions}\\
\kappa(\alpha)x^{2\alpha+1}%
\end{array}
& \longrightarrow &
\begin{array}
[c]{c}%
\text{Bessel-type functions}\\
\kappa(\alpha)x^{2\alpha+1}+M\delta(x)
\end{array}
\end{array}
\right.  \label{eq3.2}%
\end{equation}
The symbol entry (here $k$ and $\kappa$ are positive numbers depending only on
the parameters $\alpha$ and $\beta$) under each special function indicates a
non-negative (generalised) \textquotedblleft weight\textquotedblright, on the
interval $(-1,1)$ or $(0,\infty),$ involved in:

\begin{itemize}
\item[$(a)$] the orthogonality property of the special functions

\item[$(b)$] the weight coefficient in the associated differential equations.
\end{itemize}

It is important to note in this diagram that:

\begin{itemize}
\item[$(i)$] a horizontal arrow $\longrightarrow$ indicates a definition
process either by a linear combination of special functions of the same type
but of different orders, or by a linear-differential combination of special
functions of the same type and order (alternatively by an application of the
Darboux transform, see \cite{GH})

\item[$(ii)$] a vertical arrow $\downarrow$ indicates a confluent limit
process of one special function to give another special function

\item[$(iii)$] the use of the symbol $M\delta(\cdot)$ is a notational device
to indicate that the monotonic function on the real line $\mathbb{R}$ defining
the weight has a jump at an end-point of the interval concerned, of magnitude
$M>0$

\item[$(iv)$] the combination of any vertical arrow $\downarrow$ with a
horizontal arrow $\longrightarrow$ must give a consistent single entry.
\end{itemize}

Information about the Jacobi-type and Laguerre-type orthogonal polynomials,
and their associated differential equations, is given in the Everitt and
Littlejohn survey paper \cite{EL1}; see in particular the references in this
paper to the introduction of the fourth-order Laguerre-type differential
equation by H.L. and A.M. Krall, Koornwinder and by Littlejohn. The general
Laguerre-type differential equation is introduced in the paper \cite{KK} by
Koekoek and Koekoek; the order of this linear differential equation is
determined by $4+2\alpha$ with $\alpha\in\mathbb{N}_{0}=\{0,1,2,\cdots\}.$

It is significant that the general order Bessel-type functions also satisfy a
linear differential equation of order $4+2\alpha$ (with $\alpha\in
\mathbb{N}_{0}$), being an inheritance from the order of the general
Laguerre-type equation.

The purpose of this survey paper is to discuss the properties of the
Bessel-type linear differential equation in the special case when $\alpha=0,$
as given in the bottom right-hand corner of the diagram; this is the
fourth-order differential equation (\ref{eq3.1}) and involves the weight
coefficient $\kappa(0)x$; its solutions should, in some sense, have
orthogonality properties with respect to the generalised weight function
$\kappa(0)x+M\delta(0),$ where $M>0$ is the parameter appearing in the
differential equation (\ref{eq3.1}); see \cite[Section 4]{EM}.

Our knowledge of the special function solutions of the Bessel-type
differential equation (\ref{eq3.1}) is now more complete than at the time the
paper \cite{EM} was written. However, the results in \cite[Section 1,
(1.8a)]{EM}, with $\alpha=0,$ show that the function defined by%
\begin{equation}
J_{\lambda}^{0,M}(x):=[1+M(\lambda/2)^{2}]J_{0}(\lambda x)-2M(\lambda
/2)^{2}(\lambda x)^{-1}J_{1}(\lambda x)\ \text{for all}\ x\in(0,\infty),
\label{eq3.3}%
\end{equation}
is a solution of the differential equation (\ref{eq3.1}), for all $\lambda
\in\mathbb{C},$ and hence for all $\Lambda\in\mathbb{C},$ and all $M>0.$ Here:

\begin{itemize}
\item[$(i)$] the parameter $M>0$

\item[$(ii)$] the parameter $\lambda\in\mathbb{C}$

\item[$(iii)$] the spectral parameter $\Lambda$ and the parameter $M,$ in the
equation (\ref{eq3.1}), and the parameters $M$ and $\lambda,$ in the
definition (\ref{eq3.3}), are connected by the relationship%
\begin{equation}
\Lambda\equiv\Lambda(\lambda,M)=\lambda^{2}(\lambda^{2}+8M^{-1})\ \text{for
all}\ \lambda\in\mathbb{C}\ \text{and all}\ M>0 \label{eq3.4}%
\end{equation}

\item[$(iv)$] $J_{0}$ and $J_{1}$ are the classical Bessel functions (of the
first kind), see \cite[Chapter III]{GNW}.
\end{itemize}

Similar arguments to the methods given in \cite{EM} show that the function
defined by%
\begin{equation}
Y_{\lambda}^{0,M}(x):=[1+M(\lambda/2)^{2}]Y_{0}(\lambda x)-2M(\lambda
/2)^{2}(\lambda x)^{-1}Y_{1}(\lambda x)\ \text{for all}\ x\in(0,\infty),
\label{eq3.5}%
\end{equation}
is also a solution of the differential equation (\ref{eq3.1}), for all
$\lambda\in\mathbb{C},$ and hence for all $\Lambda\in\mathbb{C}$ and all
$M>0;$ here, again, $Y_{0}$ and $Y_{1}$ are classical Bessel functions (of the
second kind), see \cite[Chapter III]{GNW}.

The earlier studies of the fourth-order differential equation (\ref{eq3.1})
failed to find any explicit form of two linearly independent solutions,
additional to the solutions $J_{\lambda}^{0,M}$ and $Y_{\lambda}^{0,M}.$
However, results of van Hoeij, see \cite{MvH1} and \cite{MvH}, using the
computer algebra program Maple have yielded the required two additional
solutions, here given the notation $I_{\lambda}^{0,M}$ and $K_{\lambda}%
^{0.M},$ with explicit representation in terms of the classical modified
Bessel functions $I_{0},K_{0}$ and $I_{1},K_{1}.$ These two additional
solutions are defined as follows, where as far as possible we have followed
the notation used for the solutions $J_{\lambda}^{0,M}$ and $Y_{\lambda}%
^{0,M},$

\begin{itemize}
\item[$(i)$] given $\lambda\in\mathbb{C}$, with $\arg(\lambda)\in\lbrack
0,2\pi),$ $M\in(0,\infty)$ and using the principal value of $\sqrt{\cdot},$
define%
\begin{equation}
c\equiv c(\lambda,M):=\sqrt{\lambda^{2}+8M^{-1}}\ \text{and}\ d\equiv
d(\lambda,M):=1+M(\lambda/2)^{2} \label{eq3.5a}%
\end{equation}

\item[$(ii)$] define the solution $I_{\lambda}^{0,M}$, for all $x\in
(0,\infty),$%
\begin{equation}%
\begin{array}
[c]{lll}%
I_{\lambda}^{0,M}(x) & := & -dI_{0}(cx)+\tfrac{1}{2}cMx^{-1}I_{1}(cx)\\
& := & -[1+M(\lambda/2)^{2}]I_{0}\left(  x\sqrt{\lambda^{2}+8M^{-1}}\right) \\
&  & +\sqrt{M\left(  2+M(\lambda/2)^{2}\right)  }x^{-1}I_{1}\left(
x\sqrt{\lambda^{2}+8M^{-1}}\right)
\end{array}
\label{eq3.5b}%
\end{equation}

\item[$(iii)$] define the solution $K_{\lambda}^{0,M}$, for all $x\in
(0,\infty),$%
\begin{equation}%
\begin{array}
[c]{lll}%
I_{\lambda}^{0,M}(x) & := & dK_{0}(cx)+\tfrac{1}{2}cMx^{-1}K_{1}(cx)\\
& := & [1+M(\lambda/2)^{2}]K_{0}\left(  x\sqrt{\lambda^{2}+8M^{-1}}\right) \\
&  & +\sqrt{M\left(  2+M(\lambda/2)^{2}\right)  }x^{-1}K_{1}\left(
x\sqrt{\lambda^{2}+8M^{-1}}\right)
\end{array}
\label{eq3.5c}%
\end{equation}

\end{itemize}

\begin{remark}
\label{rem3.1}We have

\begin{enumerate}
\item The four linearly independent solutions $J_{\lambda}^{0.M},Y_{\lambda
}^{0.M},I_{\lambda}^{0.M},K_{\lambda}^{0.M}$ provide a basis for all solutions
of the original differential equation (\ref{eq3.1}), subject to the
$(\Lambda,\lambda)$ connection given in (\ref{eq3.4}).

\item These four solutions are real-valued on their domain $(0,\infty)$ for
all $\lambda\in\mathbb{R}.$

\item The domain $(0,\infty)$ of the solutions $J_{\lambda}^{0,M}$ and
$I_{\lambda}^{0,M}$ can be extended to the closed half-line $[0,\infty)$ with
the properties%
\[
J_{\lambda}^{0,M}(0)=I_{\lambda}^{0,M}(0)=1\ \text{for all}\ \lambda
\in\mathbb{R}\ \text{and all}\ M\in(0,\infty).
\]

\end{enumerate}
\end{remark}

The classical Bessel differential equation, with order $\alpha=0,$ written in
a form comparable to the fourth-order equation (\ref{eq3.1}), is best taken
from the left-hand bottom corner of the diagram (\ref{eq3.2}); from
\cite[Section 1, (1.2)]{EM} with $\alpha=0$ we obtain%
\begin{equation}
-(xy^{\prime}(x))^{\prime}=\lambda^{2}xy(x)\ \text{for all}\ x\in(0,\infty);
\label{eq3.6}%
\end{equation}
here $\lambda\in\mathbb{C}$ is the spectral parameter. It is to be observed
that, formally, if the fourth-order Bessel-type equation (\ref{eq3.1}) is
multiplied by the parameter $M>0$ and then $M$ tends to zero, we obtain
essentially the classical Bessel equation of order zero (\ref{eq3.6}), on
using the spectral relationship (\ref{eq3.4}) between the parameters $\lambda$
and $\Lambda.$ This Bessel differential equation (\ref{eq3.6}) has solutions
$J_{1}(\lambda x)$ and $Y_{1}(\lambda x)$ for all $x\in(0,\infty)$ and all
$\lambda\in\mathbb{C}.$

For the need to apply the Frobenius series method of solution we also consider
the differential equation (\ref{eq3.1}) on the complex plane $\mathbb{C}$:%
\begin{equation}%
\begin{array}
[c]{r}%
w^{(4)}(z)+2z^{-1}w^{(3)}(z)-(9z^{-2}+8M^{-1})w^{\prime\prime}(z)\\
+(9z^{-3}-8M^{-1}z^{-1})w^{\prime}(z)-\Lambda w(z)=0
\end{array}
\label{eq3.7}%
\end{equation}
for all $z\in\mathbb{C}.$ In this form the equation has a \textit{regular
singularity} at the origin $0,$ and an \textit{irregular singularity} at the
point at infinity $\infty$ of the complex plane $\mathbb{C};$ all other points
of the plane are \textit{regular} or \textit{ordinary} points for the
differential equation. It should be noted that the classical Bessel
differential equation (\ref{eq3.6}) has the same classification when
considered in the complex plane $\mathbb{C}.$

A calculation shows that the Frobenius indicial roots for the regular
singularity of the differential equation (\ref{eq3.7}) at the origin $0,$ are
$\{4,2,0,-2\}.$ The application of the Frobenius series method, using the
computer programs \cite{BELM} and Maple (see \cite{MvH}), yield four linearly
independent series solutions of (\ref{eq3.7}), each with infinite radius of
convergence in the complex plane $\mathbb{C}.$ If these solutions are labelled
to hold for the Bessel-type differential equation (\ref{eq3.1}) then we have
four solutions $\{y_{r}(\cdot,\Lambda,M):r=4,2,0,-2\},$ to accord with the
indicial roots, to give the theorem, see \cite[Section 3]{DEHLM}:

\begin{theorem}
\label{th3.1}For all $\Lambda\in\mathbb{C}$ and all $M>0,$ the differential
equation $(\ref{eq3.1})$ has four linearly independent solutions
$\{y_{r}(\cdot,\Lambda,M):r=4,2,0,-2\}$, defined on $(0,\infty)\times
\mathbb{C},$ with the following series properties as $x\rightarrow0^{+},$
where the $O$-terms depend upon the complex spectral parameter $\Lambda$ and
the parameter $M,$%
\begin{equation}
\left\{
\begin{array}
[c]{lll}%
y_{4}(x,\Lambda,M) & = & x^{4}+\tfrac{1}{3}M^{-1}x^{6}+O(x^{8})\\
y_{2}(x,\Lambda,M) & = & kx^{2}+O(x^{4}\left\vert \ln(x)\right\vert )\\
y_{0}(x,\Lambda,M) & = & l+O(x^{4}\left\vert \ln(x)\right\vert )\\
y_{-2}(x,\Lambda,M) & = & mx^{-2}+O(\left\vert \ln(x)\right\vert ).
\end{array}
\right.  \label{eq3.8}%
\end{equation}
Here the fixed numbers $k,l,m\in\mathbb{R}$ and are independent of the
parameters $\Lambda$ and $M;$ these numbers are produced by the Frobenius
computer program \cite{BELM} and have the explicit values:%
\[
k=-(27720)^{-1},\ l=(174636000)^{-1},\ m=-(9779616000)^{-1}.
\]

\end{theorem}

\section{Higher-order differential equations\label{sec3a}}

As mentioned in Section \ref{sec3} above there exist Bessel-type linear
differential equations of all even-orders $4+2\alpha$, where $\alpha
\in\mathbb{N}_{0}$ is any non-negative integer. The definition and some
properties of these differential equations, and the associated Bessel-type
functions, are considered in detail in \cite[Sections 2 and 3]{EM}.

Here we give the form of the sixth-order and eighth-order differential
equations, as given in \cite[Section 1, (1.10b) and (1.10c).]{EM}. (Note that
there is a printing error in the display (1.10b); the numerical factor $255$
is to be replaced by 225. Also printing errors in the display (1.10c) which
are now to be corrected using the form of the differential equation
(\ref{eq3.8b}) below.)

\begin{itemize}
\item[$(i)$] The sixth-order equation derived from the corrected differential
expression for \cite[Section 1, (1.9) and (1.10b)]{EM} is%
\begin{equation}%
\begin{array}
[c]{r}%
-(x^{3}y^{(3)}(x))^{(3)}+(33xy^{\prime\prime}(x))^{\prime\prime}%
-((225x^{-1}-96M^{-1}x^{3})y^{\prime}(x))^{\prime}\\
=(\lambda^{6}+M^{-1}2^{4}\left(  3!\right)  \lambda^{2})x^{3}y(x)\ \text{for
all}\ x\in(0,\infty),
\end{array}
\label{eq3.8a}%
\end{equation}
where, as before, the parameters $M\in(0,\infty)$ and $\lambda\in\mathbb{C}.$

When this equation is considered in the complex plane $\mathbb{C}$ the
Frobenius indicial roots for the regular singularity at the origin $0$ are
$\{6,4,2,0,-2,-4\},$ using the methods provided by \cite{MvH}.

\item[$(ii)$] The eighth-order equation derived from the corrected
differential expression for \cite[Section 1, (1.9) and (1.10c)]{EM} is
\begin{equation}%
\begin{array}
[c]{r}%
(x^{5}y^{(4)}(x))^{(4)}-(78x^{3}y^{(3)}(x))^{(3)}+(1809xy^{\prime\prime
}(x))^{\prime\prime}\\
-((11025x^{-1}-2^{6}(4!)M^{-1}x^{5})y^{\prime}(x))^{\prime}\\
=(\lambda^{8}+M^{-1}2^{6}(4!)\lambda^{2})x^{5}y(x)\ \text{for all}%
\ x\in(0,\infty),
\end{array}
\label{eq3.8b}%
\end{equation}
where, as before, the parameters $M\in(0,\infty)$ and $\lambda\in\mathbb{C}.$

When this equation is considered in the complex plane $\mathbb{C}$ the
Frobenius indicial roots for the regular singularity at the origin $0$ are
$\{8,6,4,2,0,-2,-4,-6\},$ using the methods provided by \cite{MvH}.
\end{itemize}

We note that, formally, if the equations in (\ref{eq3.8a}) and (\ref{eq3.8b})
are multiplied by $M>0,$ and then letting $M$ tend to zero we obtain,
respectively, the two Sturm-Liouville differential equations, see
\cite[Section 1, (1.2)]{EM},%
\begin{equation}
-(x^{3}y^{\prime}(x))^{\prime}=\lambda^{2}x^{3}y(x)\ \text{and}\ -(x^{5}%
y^{\prime}(x))^{\prime}=\lambda^{2}x^{5}y(x)\ \text{for all}\ x\in(0,\infty).
\label{eq3.8c}%
\end{equation}
For the solutions of these equations in classical Bessel functions see
\cite[Section 1, (1.4)]{EM}.

\section{The fourth-order differential expression $L_{M}$\label{sec4}}

We define the differential expression $L_{M}$ with domain $D(L_{M})$ as
follows:%
\begin{equation}
D(L_{M}):=\{f:(0,\infty)\rightarrow\mathbb{C}:f^{(r)}\in AC_{\text{loc}%
}(0,\infty)\;\text{for}\;r=0,1,2,3\}, \label{eq4.1}%
\end{equation}
and for all $f\in D(L_{M})$%
\begin{equation}
L_{M}[f](x):=(xf^{\prime\prime}(x))^{\prime\prime}-\left(  (9x^{-1}%
+8M^{-1}x)f^{\prime}(x)\right)  ^{\prime}\;(x\in(0,\infty)); \label{eq4.2}%
\end{equation}
it follows that%
\begin{equation}
L_{M}:D(L_{M})\rightarrow L_{\text{loc}}^{1}(0,\infty). \label{eq4.3}%
\end{equation}

The Green's formula for $L_{M}$ on any compact interval $[\alpha,\beta
]\subset(0,+\infty)$ is given by%
\begin{equation}
\int_{\alpha}^{\beta}\left\{  \overline{g}(x)L_{M}[f](x)-f(x)\overline
{L_{M}[g]}(x)\right\}  dx=\left.  [f,g](x)\right\vert _{\alpha}^{\beta},
\label{eq4.5}%
\end{equation}
where the symplectic form $[\cdot,\cdot](\cdot):D(L_{M})\times D(L_{M}%
)\times(0,+\infty)\rightarrow\mathbb{C}$ is defined by%
\begin{align}
\lbrack f,g](x)  &  :=\overline{g}(x)(xf^{\prime\prime}(x))^{\prime
}-(x\overline{g}^{\prime\prime}(x))^{\prime}f(x)\nonumber\\
&  \quad-x\left(  \overline{g}^{\prime}(x)f^{\prime\prime}(x)-\overline
{g}^{\prime\prime}(x)f^{\prime}(x)\right) \nonumber\\
&  \quad-\left(  9x^{-1}+8M^{-1}x\right)  \left(  \overline{g}(x)f^{\prime
}(x)-\overline{g}^{\prime}(x)f(x)\right)  . \label{eq4.6}%
\end{align}

The Dirichlet formula for $L_{M}$ on any compact interval $[\alpha
,\beta]\subset(0,+\infty)$ is given by%
\begin{equation}%
\begin{array}
[c]{c}%
{\displaystyle\int_{\alpha}^{\beta}}
\left\{  xf^{\prime\prime}(x)\overline{g}^{\prime\prime}(x)+\left(
9x^{-1}+8M^{-1}x\right)  f^{\prime}(x)\overline{g}^{\prime}(x)\right\}  dx\\
=\left.  [f,g]_{D}(x)\right\vert _{\alpha}^{\beta}+%
{\displaystyle\int_{\alpha}^{\beta}}
L_{M}[f](x)\overline{g}(x)\,dx,
\end{array}
\label{eq4.7}%
\end{equation}
where the Dirichlet form $[\cdot,\cdot]_{D}:D(L_{M})\times D_{0}(L_{M}%
)\times(0,+\infty)\rightarrow\mathbb{C}$ is defined by, for $f\in D(L_{M})$
and $g\in D_{0}(L_{M}),$ with%
\begin{equation}
D_{0}(L_{M}):=\{g:(0,+\infty)\rightarrow\mathbb{C}:g^{(r)}\in AC_{\text{loc}%
}(0,+\infty)\;\text{for}\;r=0,1\} \label{eq4.8}%
\end{equation}
and%
\begin{equation}
\lbrack f,g]_{D}(x):=-\overline{g}(x)\left(  xf^{\prime\prime}(x)\right)
^{\prime}+\overline{g}^{\prime}(x)xf^{\prime\prime}(x)+\overline{g}(x)\left(
9x^{-1}+8M^{-1}x\right)  f^{\prime}(x). \label{eq4.9}%
\end{equation}

\section{Hilbert function spaces\label{sec5}}

The spectral properties of the fourth-order Bessel differential equation%
\begin{equation}
(xy^{\prime\prime}(x))^{\prime\prime}-((9x^{-1}+8M^{-1}x)y^{\prime
}(x))^{\prime}=\Lambda xy(x)\ \text{for all}\ x\in(0,\infty), \label{eq5.1}%
\end{equation}
with $\Lambda\in\mathbb{C}$ as the spectral parameter, are considered in two
Hilbert function spaces:

\begin{enumerate}
\item The Lebesgue weighted space%
\begin{equation}
L^{2}((0,\infty);x):=\left\{  f:(0,+\infty)\rightarrow\mathbb{C}:\int
_{0}^{\infty}x\left\vert f(x)\right\vert ^{2}\,dx<+\infty\right\}
\label{eq5.2}%
\end{equation}
with inner-product and norm defined by, for all $f,g\in L^{2}((0,\infty);x),$%
\begin{equation}
(f,g):=\int_{0}^{\infty}xf(x)\overline{g}(x)~\,dx\quad\text{and}%
\quad\left\Vert f\right\Vert :=(f,f)^{1/2}. \label{eq5.3}%
\end{equation}
This space takes into account the weight function $x$ on the right-hand side
of (\ref{eq5.1}).

\item The Lebesgue-Stieltjes jump space $L^{2}([0,\infty);m_{k})$, as
suggested by the results in \cite[Section 4]{EM}.

Let the monotonic non-decreasing function $\hat{m}_{k}:[0,\infty
)\rightarrow\lbrack0,\infty)$ be defined by, where $k>0$ is a real parameter,%
\begin{align*}
\hat{m}_{k}(x)  &  =-k\ \text{for }x=0\\
&  =x^{2}/2\ \text{for all}\ x\in(0,+\infty).
\end{align*}
Then $\hat{m}_{k}$ generates a Baire measure $m_{k}$ on the $\sigma$-algebra
$\mathcal{B}$ of Borel sets on the interval $[0,\infty)$; in turn this measure
generates a Lebesgue-Stieltjes integral for Borel measurable functions.

The Hilbert function space $L^{2}([0,\infty);m_{k})$ is defined on all
functions with the properties:%
\[%
\begin{array}
[c]{ll}%
(i) & f:[0,\infty)\rightarrow\mathbb{C\ }\text{and is Borel measurable
on}\ [0,\infty)\\
(ii) &
{\displaystyle\int_{0}^{\infty}}
x\left\vert f(x)\right\vert ^{2}dx<+\infty.
\end{array}
\]
The norm and inner-product in $L^{2}([0,\infty);m_{k})$ are defined by%
\begin{equation}
\left\Vert f\right\Vert _{k}^{2}:=\int_{[0,\infty)}\left\vert f(x)\right\vert
^{2}dm_{k}(x)=k\left\vert f(0)\right\vert ^{2}+\int_{0}^{\infty}x\left\vert
f(x)\right\vert ^{2}dx \label{eq5.4}%
\end{equation}
and%
\begin{equation}
(f,g)_{k}:=\int_{[0,\infty)}f(x)\overline{g}(x)~dm_{k}(x)=kf(0)\overline
{g}(0)+\int_{0}^{\infty}xf(x)\overline{g}(x)~dx. \label{eq5.5}%
\end{equation}
Note that the first integrals in both these definitions are Lebesgue-Stieltjes
integrals taken over the set $[0,\infty),$ whilst the second integrals can be
taken as Lebesgue integrals.
\end{enumerate}

\section{Differential operators generated by $L_{M}$\label{sec5a}}

The Lagrange symmetric differential expression $L_{M}$ generates self-adjoint
operators in both the Hilbert function spaces $L^{2}((0,\infty);x)$, and in
$L^{2}([0,\infty);m_{k})$ for all $k\in(0,\infty).$

In the space $L^{2}((0,\infty);x)$ the expression $L_{M}$ generates a
continuum $\{T\}$ of self-adjoint operators, including the significant
Friedrichs operator $F$; these properties are developed and considered in
Sections \ref{sec6} to \ref{sec11} below.

For each $k\in(0,\infty)$ the expression $L_{M}$ generates a unique
self-adjoint operator $S_{k}$ in the space $L^{2}([0,\infty);m_{k});$ the
properties of this operator are considered in Sections \ref{sec12} and
\ref{sec13}.

\section{Differential operators in $L^{2}((0,\infty);x)$\label{sec6}}

The maximal and the minimal differential operators, denoted respectively
$T_{1}$ and $T_{0},$ as generated by the differential expression $L_{M}$ in
the Hilbert function space $L^{2}((0,\infty);x),$ are defined as follows, see
\cite[Chapter V, Section 17]{MAN}:

\begin{itemize}
\item[$(i)$] $T_{1}:D(T_{1})\subset L^{2}((0,\infty);x)\rightarrow
L^{2}((0,\infty);x)$ by%
\begin{equation}
D(T_{1}):=\{f\in D(L_{M}):f,x^{-1}L_{M}(f)\in L^{2}((0,\infty);x)\}
\label{eq6.1}%
\end{equation}
and%
\begin{equation}
T_{1}f:=x^{-1}L_{M}(f)\ \text{for all}\ f\in D(T_{1}). \label{eq6.2}%
\end{equation}
From the Green's formula (\ref{eq4.5}) if follows that the limits%
\begin{equation}
\lbrack f,g](0^{+}):=\lim_{x\rightarrow0}[f,g](x)\ \text{and}\ [f,g](\infty
):=\lim_{x\rightarrow\infty}[f,g](x) \label{eq6.3}%
\end{equation}
both exist and are finite in $\mathbb{C}$ for all $f,g\in D(T_{1}).$

\item[$(ii)$] $T_{0}:D(T_{0})\subset L^{2}((0,\infty);x)\rightarrow
L^{2}((0,\infty);x)$ by%
\begin{equation}%
\begin{array}
[c]{r}%
D(T_{0}):=\{f\in D(T_{1}):\lim_{x\rightarrow0}[f,g](x)=0\ \text{and\qquad
\qquad}\\
\lim_{x\rightarrow\infty}[f,g](x)=0\text{ for all}\ f,g\in D(T_{1})\},
\end{array}
\label{eq6.4}%
\end{equation}

and%
\begin{equation}
T_{0}f:=x^{-1}L_{M}(f)\ \text{for all}\ f\in D(T_{0}). \label{eq6.5}%
\end{equation}

\end{itemize}

From standard results we have the operator properties, see \cite[Chapter
V]{MAN},%
\begin{equation}
T_{0}\subseteq T_{1},T_{0}^{\ast}=T_{1}\ \text{and}\ T_{1}^{\ast}=T_{0},
\label{eq6.6}%
\end{equation}
thereby noting that both $T_{0}$ and $T_{1}$ are closed linear operators in
$L^{2}((0,\infty);x).$

\section{Self-adjoint operators in $L^{2}((0,\infty);x)$\label{sec7}}

In the weighted space $L^{2}((0,\infty);x)$ the Lagrange symmetric (formally
self-adjoint) differential expression has the following endpoint
classifications at the singular endpoints $0$ and $+\infty$ (for additional
details see \cite[Section 6]{DEHLM}):

\begin{itemize}
\item[$(i)$] At $0^{+}$ the singular endpoint is limit-3 in $L^{2}%
((0,\infty);x)$

\item[$(ii)$] At $+\infty$ the singular endpoint is Dirichlet and strong
limit-2 in $L^{2}((0,\infty);x).$
\end{itemize}

Based on this information the self-adjoint extensions of the closed symmetric
operator $T_{0}$ are determined by the GKN-theorem on singular boundary
conditions as given in \cite[Chapter V]{MAN} and \cite{EM1}. In particular,
for the operators $T_{0}$ and $T_{1},$ any self-adjoint operator $T=T^{\ast}$
generated by $L_{M}$ in $L^{2}((0,\infty);x)$ is a one-dimensional extension
of $T_{0}$ or, equivalently, a one-dimensonal restriction of $T_{1}.$ Let the
domain $D(T)$ as a restriction of the domain $D(T_{1})$ be determined by%
\begin{equation}
D(T):=\{f\in D(T_{1}):[f,\varphi](0)=0\}, \label{eq7.1}%
\end{equation}
where the function $\varphi\in D(T_{1})$ is a non-null element of the quotient
space $D(T_{1})\diagup D(T_{0})$ which satisfies the GKN symmetry condition%
\begin{equation}
\lbrack\varphi,\varphi](0)=0. \label{eq7.2}%
\end{equation}
Then the differential operator $T$ defined by%
\begin{equation}
Tf:=x^{-1}L_{M}[f]\ \text{for all}\ f\in D(T) \label{eq7.3}%
\end{equation}
satisfies $T^{\ast}=T,$ and is self-adjoint in the Hilbert space
$L^{2}((0,\infty);x).$ All such self-adjoint operators are determined in this
way on making an appropriate choice of the boundary condition function
$\varphi.$

\section{Boundary properties at $0^{+}$\label{sec8}}

The results of the following theorem are essential to obtaining the explicit
forms of the boundary conditions at $0^{+}$ to determine all self-adjoint
extensions of $T_{0}.$

\begin{theorem}
\label{th8.1}Let $f\in D(T_{1});$ then the values of $f,f^{\prime}%
,f^{\prime\prime}$ can be defined at the point $0$ so that the following
results hold:

\begin{itemize}
\item[$(i)$] $f\in AC[0,1]$

\item[$(ii)$] $f^{\prime}\in AC[0,1]$ and $f^{\prime}(0)=0$

\item[$(iii)$] $f^{\prime\prime}\in AC_{\text{loc}}(0,1]$ and $f^{\prime
\prime}\in C[0,1]$

\item[$(iv)$] $f^{(3)}\in AC_{\text{loc}}(0,1]$ and $\lim_{x\rightarrow0^{+}%
}(xf^{(3)}(x))=0.$
\end{itemize}
\end{theorem}

For the proof of this theorem see \cite[Section 8]{DEHLM}.

We consider the functions $1,x,x^{2}$ on the interval $[0,1]$ but
\textquotedblleft patched\textquotedblright, see the Naimark patching lemma
\cite[Chapter V, Section 17.3, Lemma 2]{MAN}, to zero on $[2,\infty)$ in such
a manner that the patched functions belong to the domain $D(T_{1})$; we
continue to use the symbols $1,x,x^{2}$ for the patched functions.

A calculation shows that the results given in the next lemma are satisfied:

\begin{lemma}
\label{lem8.1}The patched functions $1,x,x^{2}$ have the following limit
properties in respect of the symplectic form $[\cdot,\cdot]$ and the maximal
domain $D(T_{1})$:

\begin{itemize}
\item[$(i)$] $1,x^{2}\in D(T_{1})$ but $x\notin D(T_{1})$

\item[$(ii)$] $[1,1](0^{+})=[x,x](0^{+})=[x^{2},x^{2}](0^{+})=0$

\item[$(iii)$] $[x,x^{2}](0^{+})=0$ and $[1,x^{2}](0^{+})=16$

\item[$(iv)$] $[1,x](0^{+})$ does not exist.
\end{itemize}
\end{lemma}

The lemmas and corollaries now given below are taken from \cite[Section
9]{DEHLM}, where proofs are given in detail.

The results of Theorem \ref{th8.1} and Lemma \ref{lem8.1} now provide a basis
for the two-dimensional quotient space $D(T_{1})/D(T_{0})$;%
\begin{equation}
D(T_{1})/D(T_{0})=\mathrm{span}\{1,x^{2}\}=\{a+bx^{2}:a,b\in\mathbb{C}\}.
\label{eq8.1}%
\end{equation}
The linear independence of the functions $\{1,x^{2}\}$ within the the quotient
space follows from the property $[1,x^{2}](0^{+})=16\neq0.$

A calculation now gives, recall Theorem \ref{th8.1},

\begin{lemma}
\label{lem8.2}Let $f\in D(T_{1});$ then the following identities hold:

\begin{itemize}
\item[$(i)$] $[f,1](0^{+})=-8f^{\prime\prime}(0)$

\item[$(ii)$] $[f,x^{2}](0^{+})=16f(0).$
\end{itemize}
\end{lemma}

Similarly we have

\begin{lemma}
\label{lem8.3}Let $f,g\in D(T_{1});$ then

\begin{itemize}
\item[$(i)$] $[f,g](0^{+})=8\left[  f(0)\overline{g}^{\prime\prime
}(0)-f^{\prime\prime}(0)\overline{g}(0)\right]  $

\item[$(ii)$] $\left[  f,g\right]  _{D}(0^{+})=8f^{\prime\prime}%
(0)\overline{g}(0).$
\end{itemize}
\end{lemma}

We have the corollaries:

\begin{corollary}
\label{cor8.1}The domain of the minimal operator $T_{0}$ is determined
explicitly by%
\begin{equation}
D(T_{0})=\{f\in D(T_{1}):f(0)=0\ \text{and}\ f^{\prime\prime}(0)=0\}.
\label{eq8.2}%
\end{equation}

\end{corollary}

\begin{corollary}
\label{cor8.3}For all $f\in D(T_{1})$
\begin{equation}
\int_{0}^{\infty}\left\{  x\left\vert f^{\prime\prime}(x)\right\vert
^{2}+\left(  9x^{-1}+8M^{-1}x\right)  \left\vert f^{\prime}(x)\right\vert
^{2}\right\}  dx<\infty. \label{eq8.4}%
\end{equation}

\end{corollary}

\begin{corollary}
\label{cor8.4}For all $f,g\in D(T_{1})$ the Dirichlet formula takes the form%
\begin{equation}
(T_{1}f,g)=8f^{\prime\prime}(0)\overline{g}(0)+\int_{0}^{\infty}\left\{
xf^{\prime\prime}(x)\overline{g}^{\prime\prime}(x)+\left(  9x^{-1}%
+8M^{-1}x\right)  f^{\prime}(x)\overline{g}^{\prime}(x)\right\}  dx.
\label{eq8.5}%
\end{equation}

\end{corollary}

\section{Explicit boundary condition functions at $0^{+}$\label{sec9}}

We can now determine all forms of the boundary condition function $\varphi$
satisfying the symmetry condition (\ref{eq7.2}) to determine the domain of all
self-adjoint extensions $T$ of the minimal operator $T_{0}.$

\begin{lemma}
\label{lem9.1}All self-adjoint extensions $T$ of $T_{0}$ generated by the
differential expression $L_{M}$ in $L^{2}((0,\infty;x)$ are determined by,
using the patched functions $1,x^{2},$%
\begin{equation}
D(T):=\{f\in D(T_{1}):[f,\varphi](0^{+})=0\ \text{where}\
\begin{array}
[t]{ll}%
(i) & \varphi(x)=\alpha+\beta x^{2}\\
(ii) & \alpha,\beta\in\mathbb{R}\ \text{and}\ \alpha^{2}+\beta^{2}\neq0\}.
\end{array}
\ \label{eq9.1}%
\end{equation}
and%
\begin{equation}
\left(  Tf\right)  (x):=x^{-1}L_{M}(f)(x)\ \text{for all}\ x\in(0,\infty
)\ \text{and all}\ f\in D(T). \label{eq9.2}%
\end{equation}

\end{lemma}

There is an equivalent form of this last result, using the results of Lemma
\ref{lem8.2}:

\begin{lemma}
\label{lem9.2}All self-adjoint extensions $T$ of $T_{0}$ generated by the
differential expression $L_{M}$ in $L^{2}((0,\infty;x)$ are determined by%
\begin{equation}
D(T):=\{f\in D(T_{1}):%
\begin{array}
[t]{ll}%
(i) & -\alpha f^{\prime\prime}(0)+2\beta f(0)=0\\
(ii) & \alpha,\beta\in\mathbb{R}\ \text{and}\ \alpha^{2}+\beta^{2}\neq0\}.
\end{array}
\label{eq9.3}%
\end{equation}
and%
\begin{equation}
\left(  Tf\right)  (x):=x^{-1}L_{M}(f)(x)\ \text{for all}\ x\in(0,\infty
)\ \text{and all}\ f\in D(T). \label{eq9.4}%
\end{equation}

\end{lemma}

\begin{remark}
\label{rem9.1}We note the two special cases:

\begin{itemize}
\item[$(i)$] When $\alpha=0$ the boundary condition is $f(0)=0;$ this boundary
condition plays a special role, and gives an explicit form of the domain of
the Friedrichs extension $F$ of $T_{0};$ see Section \ref{sec11} below.

\item[$(ii)$] When $\beta=0$ the boundary condition is $f^{\prime\prime
}(0)=0.$
\end{itemize}
\end{remark}

\section{Spectral properties of the fourth-order Bessel-type
operators\label{sec10}}

\begin{theorem}
\label{th10.1}The minimal operator $T_{0},$ defined in $(\ref{eq6.4})$ and
$(\ref{eq6.5}),$ is bounded below in the space $L^{2}((0,\infty);x)$ by the
null operator $O,$ \textit{i.e.}%
\begin{equation}
(T_{0}f,f)\geq0\ \text{for all}\ f\in D(T_{0}). \label{eq10.1}%
\end{equation}

\end{theorem}

\begin{proof}
Since $T_{0}$ is a restriction of the maximal operator $T_{1}$ the result of
Corollary \ref{cor8.4} can be applied to give, using also Corollary
\ref{cor8.1},%
\[%
\begin{array}
[c]{c}%
(T_{0}f,f)=8f^{\prime\prime}(0)\overline{f}(0)+%
{\displaystyle\int_{0}^{\infty}}
\left\{  x\left\vert f^{\prime\prime}(x)\right\vert ^{2}+\left(
9x^{-1}+8M^{-1}x\right)  \left\vert f^{\prime}(x)\right\vert ^{2}\right\}
dx\\
=%
{\displaystyle\int_{0}^{\infty}}
\left\{  x\left\vert f^{\prime\prime}(x)\right\vert ^{2}+\left(
9x^{-1}+8M^{-1}x\right)  \left\vert f^{\prime}(x)\right\vert ^{2}\right\}
dx\geq0\
\end{array}
\]
for all $f\in D(T_{0}).$
\end{proof}

\begin{theorem}
\label{th10.2}

\begin{enumerate}
\item Let $T$ be a self-adjoint extension of $T_{0};$ then:

\begin{itemize}
\item[$(i)$] The essential spectrum $\sigma_{\text{ess}}(T)$ is given by%
\begin{equation}
\sigma_{\text{ess}}(T)=\sigma_{\text{cont}}(T)=[0,\infty). \label{eq10.2}%
\end{equation}

\item[$(ii)$] There are no embedded eigenvalues of $T$ in the essential spectrum.

\item[$(iii)$] $T$ has at most one eigenvalue; if this eigenvalue is present
then it is simple and lies in the interval $(-\infty,0).$
\end{itemize}

\item Every point $\mu\in(-\infty,0)$ is the eigenvalue of some unique
self-adjoint extension $T$ of $T_{0}.$
\end{enumerate}
\end{theorem}

\begin{proof}
The proof of this theorem is given in detail in \cite[Section 13]{DEHLM}.
\end{proof}

\section{The Friedrichs extension $F$\label{sec11}}

The closed symmetric operator $T_{0}$ is bounded below in $L^{2}%
((0,\infty);x),$ see Theorem \ref{th10.1}, and the general theory of such
operators implies the existence of a distinguished self-adjoint extension $F,$
called the Friedrichs extension of $T_{0}.$

This Friedrichs operator has the properties:

\begin{itemize}
\item[$(i)$] $T_{0}\subset F=F^{\ast}\subset T_{1}$

\item[$(ii)$] $D(F)=\{f\in D(T_{1}):f(0)=0\}$

\item[$(iii)$] The essential spectrum $\sigma_{\text{ess}}(F)$ is given by%
\begin{equation}
\sigma_{\text{ess}}(F)=\sigma_{\text{cont}}(F)=[0,\infty) \label{eq11.1}%
\end{equation}

\item[$(iv)$] $F$ has no eigenvalues.
\end{itemize}

For a discussion of the definition and properties of this Friedrichs extension
see \cite[Section 15]{DEHLM}.

\section{Self-adjoint operator $S_{k}$ in $L^{2}([0,\infty);m_{k}%
)$\label{sec12}}

In this section, given any $k\in(0,\infty),$ we define the operator $S_{k}$
generated by the differential expression $L_{M}$ in the Hilbert function space
$L^{2}([0,\infty);m_{k}),$ where this space is defined in Section \ref{sec5} above.

\begin{definition}
\label{def12.1}Let $k\in(0,\infty)$ be given$;$ then the operator $S_{k}$%
\begin{equation}
S_{k}:D(S_{k})\subset L^{2}([0,\infty);m_{k})\rightarrow L^{2}([0,\infty
);m_{k}) \label{eq12.1}%
\end{equation}
is defined by $($see $(\ref{eq6.1})$ and $(\ref{eq6.2}),$ and Theorem
$\ref{th8.1}$ for the definition and properties of the domain $D(T_{1})\subset
L^{2}((0,\infty);x))$

\begin{itemize}
\item[$(i)$] $D(S_{k}):=D(T_{1})$

\item[$(ii)$] for all $f\in D(S_{k})$%
\begin{equation}
\left\{
\begin{array}
[c]{lll}%
\left(  S_{k}f\right)  (x) & := & -8k^{-1}f^{\prime\prime}(0)\ \text{for}%
\ x=0\\
& := & x^{-1}L_{M}[f](x)\ \text{for all}\ x\in(0,\infty).
\end{array}
\right.  \label{eq12.2}%
\end{equation}

\end{itemize}
\end{definition}

\begin{theorem}
\label{th12.1}For all $k\in(0,\infty)$:

\begin{itemize}
\item[$(i)$] The linear manifold $D(S_{k})$ is dense in $L^{2}([0,\infty
);m_{k}).$

\item[$(ii)$] The operator $S_{k}$ is hermitian in $L^{2}([0,\infty);m_{k}).$

\item[$(iii)$] The operator $S_{k}$ is symmetric in $L^{2}([0,\infty);m_{k}).$

\item[$(iv)$] The operator $S_{k}$ is bounded below in $L^{2}([0,\infty
);m_{k})$%
\begin{equation}
(S_{k}f,f)_{k}\geq0\ \text{for all}\ f\in D(S_{k}). \label{eq12.3}%
\end{equation}

\end{itemize}
\end{theorem}

For the proof of this theorem see \cite[Theorem 5.2]{EKLM}.

\begin{theorem}
\label{th12.2}Let $k\in(0,\infty)$ be given; then the symmetric operator
$S_{k}$ on the domain $D(S_{k})$ is self-adjoint in the Hilbert function space
$L^{2}([0,\infty);m_{k}).$
\end{theorem}

For the proof of this theorem see \cite[Theorem 5.4]{EKLM}.

\begin{theorem}
\label{th12.3}Let $k\in(0,\infty)$ be given; then the operator $S_{k}$ on the
domain $D(S_{k})$ is the unique self-adjoint operator generated by the
differential expression $L_{M}$ in the Hilbert function space $L^{2}%
([0,\infty);m_{k}).$
\end{theorem}

For the proof of this theorem see \cite[Theorem 5.5]{EKLM}.

\section{Spectral properties of the self-adjoint operator $S_{k}$%
\label{sec13}}

The spectral properties of the self-adjoint operator $S_{k}$ in $L^{2}%
([0,\infty);m_{k})$ are given by

\begin{theorem}
\label{th13.1}For any $k\in(0,\infty)$ let the self-adjoint operator $S_{k}$
in $L^{2}([0,\infty);m_{k})$ be defined as in Definition $\ref{def12.1}$
above; then the spectrum $\sigma(S_{k})$ of $S_{k}$ has the following
properties$:$

\begin{itemize}
\item[$(i)$] $S_{k}$ has no eigenvalues

\item[$(ii)$] the essential spectrum of $S_{k}$ is given by%
\begin{equation}
\sigma_{\text{ess}}(S_{k})=\sigma_{\text{cont}}(S_{k})=[0,\infty).
\label{eq13.1}%
\end{equation}

\end{itemize}
\end{theorem}

For the proof of this theorem see \cite[Theorem 6.1]{EKLM}.

\section{Distributional orthogonality relationships\label{sec14}}

Recall that from the properties of the classical Bessel function $J_{0}$ we
have the result that $J_{0}(\cdot)$ $\notin L^{2}((0,\infty);x)$. However from
\cite[Section 1, (1.7)]{EM} we have the following distributional (Schwartzian)
orthogonal relationship for the classical Bessel function $J_{0},$ in the
space $\mathcal{D}^{\prime}$ of distributions,%
\begin{equation}
\lambda\int_{0}^{\infty}xJ_{0}(\lambda x)J_{0}(\mu x)~dx=\delta(\lambda
-\mu)\ \text{for all}\ \lambda,\mu\in(0,\infty); \label{eq14.1}%
\end{equation}
here $\delta\in\mathcal{D}^{\prime}$ is the Dirac delta distribution. This is
the generalised orthogonality property for the solutions $J_{0}$ of the
classical Bessel differential equation, of order $0,$ given by (\ref{eq3.6});
this result mirrors the spectral properties of this equation, when considered
on the half-line $(0,\infty),$ in the space $L^{2}((0,\infty);x);$ in
particular the result that every self-adjoint extension $T$ of the
corresponding minimal operator $T_{0}$ has the property $\sigma_{\text{ess}%
}(T)=[0,\infty).$

The distributional proof of (\ref{eq14.1}) is discussed in the forthcoming
paper \cite{EKLM1}, where the result is also related to the properties of
infinite integrals of Bessel functions as originated by Hankel, see
\cite[Chapter XIII]{GNW}.

As above for the Bessel function $J_{0}$ we have, from the explicit
representation (\ref{eq3.3}), the fourth-order Bessel-type function
$J_{\lambda}^{0,M}\notin L^{2}([0,\infty);m_{k})$ for all $k,M\in(0,\infty).$

To obtain a distributional orthogonality for $J_{\lambda}^{0,M}$, given any
$M>0$, it is necessary to choose a special value of the parameter $k,$
\textit{i.e.} $k=M/2.$ Then it is shown in \cite[Section 4, Corollary 4.3]{EM}
that we have the following distributional (Schwartzian) orthogonal
relationship for the fourth-order Bessel-type function $J_{\lambda}^{0,M},$ in
the space $\mathcal{D}^{\prime}$ of distributions,%
\begin{equation}%
\begin{array}
[c]{r}%
\lambda\left[  1+M(\lambda/2)^{2}\right]  ^{-2}\left\{
{\displaystyle\int_{0}^{\infty}}
xJ_{\lambda}^{0,M}(x)J_{\mu}^{0,M}(x)~dx+\tfrac{1}{2}MJ_{\lambda}%
^{0,M}(0)J_{\mu}^{0,M}(0)\right\} \\
=\delta(\lambda-\mu)\ \text{for all}\ \lambda,\mu\in(0,\infty).
\end{array}
\label{eq14.2}%
\end{equation}
The distributional proof of (\ref{eq14.2}) is discussed in the forthcoming
paper \cite{EKLM1}.

As a formal representation it follows that (\ref{eq14.2}) may be written as,
using the inner-product for the space $L^{2}\left(  [0,\infty);m_{M/2}\right)
,$%
\begin{equation}
\lambda\left[  1+M(\lambda/2)^{2}\right]  ^{-2}\left(  J_{\lambda}^{0,M}%
(\cdot),J_{\mu}^{0,M}(\cdot)\right)  _{M/2}=\delta(\lambda-\mu)~\text{for
all}\ \lambda,\mu\in(0,\infty). \label{eq14.3}%
\end{equation}

As another connection between the classical Bessel (\ref{eq3.6}) and the
fourth-order Bessel-type (\ref{eq3.1}) differential equations it is to be
noted that, formally, the orthogonality result (\ref{eq14.2}) tends to the
orthogonality result (\ref{eq14.1}), as the parameter $M$ tends to zero.

\section{The generalised Hankel transform\label{sec15}}

From the general theory of symmetric integrable-square transforms given in
\cite[Chapter VIII]{ECT} one form of the classical Hankel transform, for the
Bessel function $J_{0}$ and working in the Hilbert function space
$L^{2}((0,\infty;x),$ is:

\begin{itemize}
\item[$(i)$] Let $f\in L^{2}((0,\infty);x)$ then the Hankel transform $g\in
L^{2}((0,\infty);s)$ is given by, for $s\in(0,\infty),$%
\begin{equation}
g(s)=\int_{0}^{\infty}\xi J_{0}(s\xi)f(\xi)~d\xi\label{eq15.1}%
\end{equation}
with convergence of the integral in $L^{2}((0,\infty);s)$

\item[$(ii)$] With $g\in L^{2}((0,\infty);s)$ the inverse transform, to
recover $f,$ is given by, for $x\in(0,\infty),$%
\begin{equation}
f(x)=\int_{0}^{\infty}sJ_{0}(xs)g(s)~ds \label{eq15.2}%
\end{equation}

with convergence of the integral in $L^{2}((0,\infty);x)$

\item[$(iii)$] The Parseval relation holds between $g$ and $f$%
\begin{equation}
\int_{0}^{\infty}x\left\vert f(x)\right\vert ^{2}dx=\int_{0}^{\infty
}s\left\vert g(s)\right\vert ^{2}ds. \label{eq15.3}%
\end{equation}

\end{itemize}

There is also a direct convergence form of the Hankel transform which is best
written as, starting with $f\in L^{1}((0,\infty);x),$%
\begin{equation}
f(x)=\int_{0}^{\infty}sJ_{0}(xs)~ds\int_{0}^{\infty}\xi J_{0}(s\xi)f(\xi
)~d\xi\label{eq15.4}%
\end{equation}
with $x\in(0,\infty).$ Here the integrals are Lebesgue or limits of Lebesgue
integrals as discussed in \cite[Chapter VIII]{ECT}.

There is an equivalent generalised Hankel transform involving the fourth-order
Bessel-type function $J_{\lambda}^{0,M}(\cdot)$ and working now in the Hilbert
function space $L^{2}\left(  [0,\infty);m_{M/2}\right)  ;$ note again these
results require the unique choice of $k=M/2.$

The complete discussion of the following results for the generalised Hankel
transform are to be found in the forthcoming paper \cite{EKLM1}.

To state these results the Lebesgue-Stieltjes Hilbert function space
$L^{2}((0,\infty);n)$ is required. Let the function $\hat{n}:[0,\infty
)\rightarrow\lbrack0,\infty)$ be defined by%
\begin{equation}
\hat{n}(\lambda):=\tfrac{1}{2}\lambda^{2}\left[  1+M(\lambda/2)^{2}\right]
^{-1}\ \text{for all}\ \lambda\in\lbrack0,\infty); \label{eq15.5}%
\end{equation}
then%
\[
\hat{n}^{\prime}(\lambda)=\lambda\left[  1+M(\lambda/2)^{2}\right]  ^{-2}%
\geq0\ \text{for all}\ \lambda\in\lbrack0,\infty)
\]
so that $\hat{n}$ is monotonic increasing on $[0,\infty)$ and generates a
Baire measure on the $\sigma$-algebra $\mathcal{B}$ of Borel sets on the
interval $[0,\infty).$ The Hilbert space $L^{2}((0,\infty);n)$ is then defined
as the set of all Borel measurable complex-valued functions $f$ on
$[0,\infty)$ such that%
\[
\int_{\lbrack0,\infty)}\left\vert f(\lambda)\right\vert ^{2}dn(\lambda
)<+\infty,
\]
with norm and inner-product defined by%
\begin{equation}
\left\Vert f\right\Vert _{n}^{2}:=\int_{0}^{\infty}\left\vert f(\lambda
)\right\vert ^{2}dn(\lambda)=\int_{0}^{\infty}\left\vert f(\lambda)\right\vert
^{2}\lambda\left[  1+M(\lambda/2)^{2}\right]  ^{-2}d\lambda\label{eq15.6}%
\end{equation}%
\[
(f,g)_{n}=\int_{[0,\infty)}f(\lambda)\overline{g}(\lambda)dn(\lambda)=\int
_{0}^{\infty}f(\lambda)\overline{g}(\lambda)\lambda\left[  1+M(\lambda
/2)^{2}\right]  ^{-2}d\lambda.
\]

\begin{remark}
\label{rem15.0}This norm $\left\Vert \cdot\right\Vert _{n}$ and inner-product
$(\cdot,\cdot)_{n}$ for the space $L^{2}((0,\infty);n)$ are not to be confused
with the norm $\left\Vert \cdot\right\Vert _{k}$ and inner-product
$(\cdot,\cdot)_{k},$ introduced in Section \ref{sec5}, for the space
$L^{2}([0,\infty);m_{k}).$
\end{remark}

We note that the weight function $\lambda\longmapsto\lambda\left[
1+M(\lambda/2)^{2}\right]  ^{-2}$ in the integral in (\ref{eq15.6}) is the
factor in the distributional orthogonal relationships (\ref{eq14.2}) and
(\ref{eq14.3}).

\begin{enumerate}
\item The $L^{2}$-theory of the generalised Hankel transform is given by the
following results:

\begin{theorem}
\label{th15.1}Let $f\in L^{2}\left(  [0,\infty);m_{M/2}\right)  .$ Then there
exists exactly one function $g\in L^{2}((0,\infty);n)$ with the property that%
\begin{equation}
\int_{0}^{\infty}\left\vert g(\lambda)\right\vert ^{2}dn(\lambda
)=\int_{[0,\infty)}\left\vert f(x)\right\vert ^{2}dm_{M/2}(x); \label{eq15.7}%
\end{equation}
here $g$ is defined by, for almost all $\lambda\in(0,\infty),$%
\begin{equation}
\left(  \mathcal{F}_{M}f\right)  (\lambda):=g(\lambda)=\int_{[0,\infty
)}J_{\lambda}^{0,M}(x)f(x)dm_{M/2}(x), \label{eq15.8}%
\end{equation}
thereby defining also the generalised Hankel operator
\begin{equation}
\mathcal{F}_{M}:L^{2}\left(  [0,\infty);m_{M/2}\right)  \rightarrow
L^{2}((0,\infty);n). \label{eq15.9}%
\end{equation}

In addition $g$ satisfies%
\begin{equation}
\int_{0}^{\infty}g(\lambda)dn(\lambda)=f(0). \label{eq15.10}%
\end{equation}

\end{theorem}

\begin{remark}
\label{rem15.1}Note that the result (\ref{eq15.8}) has to be interpreted as
follows, in $(i)$ and $(ii)$:

\begin{itemize}
\item[$(i)$] $%
{\displaystyle\int_{[0,X]}}
J_{\lambda}^{0,M}(x)f(x)dm_{M/2}(x)\in L^{2}((0,\infty);n)\ $for
all$\ X\in\lbrack0,\infty)$

\item[$(ii)$]
\[
\lim_{X\rightarrow\infty}\int_{0}^{\infty}\left\vert g(\lambda)-\int
_{[0,X]}J_{\lambda}^{0,M}(x)f(x)dm_{M/2}(x)\right\vert ^{2}dn(\lambda)=0.
\]

\item[$(iii)$] Note that the Cauchy-Schwarz inequality shows that if $g\in
L^{2}((0,\infty);n)$ then $g\in L^{1}((0,\infty);n).$
\end{itemize}
\end{remark}

\begin{theorem}
\label{th15.2}Let $g\in L^{2}((0,\infty);n).$ Then there is exactly one
function $f\in L^{2}\left(  [0,\infty);m_{M/2}\right)  $ with the property
that $(\ref{eq15.7})$ is satisfied; here $f$ is defined by%
\begin{equation}
\left(  \mathcal{G}_{M}g\right)  (x):=f(x):=\left\{
\begin{array}
[c]{ll}%
{\displaystyle\int_{0}^{\infty}}
g(\lambda)dn(\lambda) & \text{if}\ x=0\\
& \\%
{\displaystyle\int_{0}^{\infty}}
J_{\lambda}^{0,M}(x)g(\lambda)dn(\lambda) & \text{for}\ x\in(0,\infty),
\end{array}
\right.  \label{eq15.11}%
\end{equation}
thereby defining also the inverse generalised Hankel operator%
\begin{equation}
\mathcal{G}_{M}:L^{2}((0,\infty);n)\rightarrow L^{2}\left(  [0,\infty
);m_{M/2}\right)  . \label{eq15.12}%
\end{equation}

\end{theorem}

\begin{remark}
\label{rem15.2}Note that the result (\ref{eq15.11}) has to be interpreted as follows:

\begin{itemize}
\item[$(i)$] $%
{\displaystyle\int_{0}^{\Lambda}}
J_{\lambda}^{0,M}(x)g(\lambda)dn(\lambda)\in L^{2}\left(  [0,\infty
);m_{M/2}\right)  \ $for all$\ \Lambda\in(0,\infty)$

\item[$(ii)$]
\[
\lim_{\Lambda\rightarrow\infty}\int_{0}^{\infty}\left\vert f(x)-%
{\displaystyle\int_{0}^{\Lambda}}
J_{\lambda}^{0,M}(x)g(\lambda)dn(\lambda)\right\vert ^{2}xdx=0.
\]

\end{itemize}
\end{remark}

\item The direct convergence of the generalised Hankel transform is given by
the following results:

\begin{theorem}
\label{th15.3}Let $\gamma\in(0,\infty).$ If $f:(0,\infty)\rightarrow
\mathbb{R}$ has the property%
\begin{equation}
x\longmapsto\sqrt{x}f(x)\in L^{1}(0,\infty) \label{eq15.13}%
\end{equation}
and is of bounded variation in a neighbourhood of $\gamma,$ then%
\begin{equation}
\tfrac{1}{2}[f(\gamma+0)+f(\gamma-0)]=\int_{0}^{\infty}J_{\lambda}%
^{0,M}(\gamma)\left(  \int_{0}^{\infty}J_{\lambda}^{0,M}(x)f(x)xdx\right)
dn(\lambda). \label{eq15.14}%
\end{equation}

Let $\mu\in(0,\infty).$ If $g:(0,\infty)\rightarrow\mathbb{R}$ has the
property%
\begin{equation}
\lambda\longmapsto\frac{\sqrt{\lambda}g(\lambda)}{1+M(\lambda/2)^{2}}\in
L^{1}(0,\infty) \label{eq15.15}%
\end{equation}

and is of bounded variation in a neighbourhood of $\mu,$ then%
\begin{equation}%
\begin{array}
[c]{l}%
\tfrac{1}{2}[g(\mu+0)+g(\mu-0)]\\
\qquad\qquad=%
{\displaystyle\int_{0}^{\infty}}
J_{\mu}^{0,M}(x)\left(
{\displaystyle\int_{0}^{\infty}}
J_{\lambda}^{0,M}(x)g(\lambda)dn(\lambda)\right)  dm_{M/2}(x).
\end{array}
\label{eq15.16}%
\end{equation}

\end{theorem}

\begin{remark}
\label{rem15.3}The integrals in Theorem \ref{th15.3} are either Lebesgue
integrals or limits of such integrals over compact intervals of $(0,\infty).$
\end{remark}

\begin{corollary}
\label{cor15.1}

\begin{itemize}
\item[$(i)$] If $\gamma\in(0,\infty)$ is a point of continuity of the function
$f$ then $(\ref{eq15.8}),~(\ref{eq15.11})$ and $(\ref{eq15.14})$ imply%
\begin{align}
(\mathcal{G}_{M}(\mathcal{F}_{M}f))(\gamma)  &  =\int_{0}^{\infty}J_{\lambda
}^{0,M}(\gamma)(\mathcal{F}_{M}f)(\lambda)dn(\lambda)\nonumber\\
&  =(M/2)f(0)\int_{0}^{\infty}J_{\lambda}^{0,M}(\gamma)dn(\lambda)+\nonumber\\
&  \;\;\;\int_{0}^{\infty}J_{\lambda}^{0,M}(\gamma)\left(  \int_{0}^{\infty
}J_{\lambda}^{0,M}(x)f(x)xdx\right)  dn(\lambda)\nonumber\\
&  =f(\gamma), \label{eq15.17}%
\end{align}
since%
\begin{equation}
\int_{0}^{\infty}J_{\lambda}^{0,M}(\eta)dn(\lambda)=0\ \text{for all}\ \eta
\in(0,\infty). \label{eq15.18}%
\end{equation}

\item[$(ii)$] If $\gamma=0$ then%
\[
(\mathcal{G}_{M}(\mathcal{F}_{M}f))(0)=(M/2)f(0)\int_{0}^{\infty}%
dn(\lambda)=f(0)
\]
since%
\[
\int_{0}^{\infty}J_{\lambda}^{0,M}(0)\left(  \int_{0}^{\infty}J_{\lambda
}^{0,M}(x)f(x)xdx\right)  dn(\lambda)=0,
\]
from $(\ref{eq15.18}),$ the use of the Fubini integral theorem and noting that
$J_{\lambda}^{0,M}(0)=1.$

\item[$(iii)$] If $\mu\in(0,\infty)$ is a point of continuity of $g$ then
$(\ref{eq15.8}),~(\ref{eq15.11})$ and $(\ref{eq15.16})$ imply%
\[
(\mathcal{F}_{M}(\mathcal{G}_{M}g))(\mu)=g(\mu).
\]

\end{itemize}
\end{corollary}
\end{enumerate}

\section{The Plum partial differential equation\label{sec16}}

The Plum equation is a fourth-order linear partial differential in the
Euclidean space $\mathbb{R}^{2}$ of two dimensions, derived from a linear
partial differential expression which is connected with the fourth-order
Bessel-type ordinary differential equation.

If the Laplacian $\nabla^{2}$ partial differential expression is written in
polar co-ordinates%
\begin{equation}
\nabla^{2}=\frac{\partial^{2}}{\partial r^{2}}+\frac{1}{r}\frac{\partial
}{\partial r}+\frac{1}{r^{2}}\frac{\partial^{2}}{\partial\theta^{2}}
\label{eq16.1}%
\end{equation}
then the Plum equation has the form, with $u=u(r,\theta),$%
\begin{equation}
\nabla^{4}u-\gamma\nabla^{2}u-\frac{4\gamma}{r^{2}}u=\Lambda u. \label{eq16.2}%
\end{equation}
Here $\gamma>0$ is determined by $\gamma=8M^{-1}$ where $M>0$ is the parameter
in the fourth-order Bessel equation (\ref{eq3.1}), and $\Lambda\in\mathbb{C}$
is a spectral parameter.

Written out the equation (\ref{eq16.2}) becomes, see \cite[Section 1, (2)]%
{L1},%
\begin{align}
&  \dfrac{\partial^{4}u}{\partial r^{4}}+\dfrac{2}{r}\dfrac{\partial^{3}%
u}{\partial r^{3}}+\left(  -\dfrac{1}{r^{2}}-\gamma\right)  \dfrac
{\partial^{2}u}{\partial r^{2}}+\left(  \dfrac{1}{r^{3}}-\dfrac{\gamma}%
{r}\right)  \dfrac{\partial u}{\partial r}+\dfrac{1}{r^{4}}\dfrac{\partial
^{4}u}{\partial\theta^{4}}\nonumber\\
&  +\dfrac{2}{r^{2}}\dfrac{\partial^{4}u}{\partial\theta^{2}\partial r^{2}%
}-\dfrac{2}{r^{3}}\dfrac{\partial^{3}u}{\partial\theta^{2}\partial r}+\left(
\dfrac{4}{r^{4}}-\dfrac{\gamma}{r^{2}}\right)  \dfrac{\partial^{2}u}%
{\partial\theta^{2}}-\dfrac{4\gamma}{r^{2}}u\nonumber\\
&  =\Lambda u. \label{eq16.3}%
\end{align}

From the results given in \cite{MP} and \cite{L1} assume that a solution for
(\ref{eq16.2}) is of the separated form%
\begin{equation}
u(r,\theta)=v(r)w(\theta) \label{eq16.4}%
\end{equation}
where $w$ is required to be a solution of the second-order Sturm-Liouville
differential equation%
\begin{equation}
-w^{\prime\prime}(\theta)=4w(\theta). \label{eq16.5}%
\end{equation}
Note that $w$ is then of the general form%
\begin{equation}
w(\theta)=A\cos(2\theta)+B\sin(2\theta) \label{eq16.6}%
\end{equation}
for, say, $\theta\in\lbrack0,\pi]$ and scalars $A,B.$

Also note that the factor $4$ in the equation (\ref{eq16.5}) is critical, and
has to be fixed, for the separation method to be effective.

Substitution of (\ref{eq16.4}) into (\ref{eq16.3}) yields, see \cite[Section
1]{L1} and \cite{MP},%
\begin{align}
&  \left(  v^{(4)}(r)+\dfrac{2}{r}v^{(3)}(r)+\left(  -\dfrac{1}{r^{2}}%
-\gamma\right)  v^{\prime\prime}(r)+\left(  \dfrac{1}{r^{3}}-\dfrac{\gamma}%
{r}\right)  v^{\prime}(r)+\dfrac{16}{r^{4}}v(r)\right)  w(\theta)\nonumber\\
&  +\left(  -\dfrac{8}{r^{2}}v^{\prime\prime}(r)+\dfrac{8}{r^{3}}v^{\prime
}(r)+\left(  -\dfrac{16}{r^{4}}+\dfrac{4\gamma}{r^{2}}\right)  v(r)-\dfrac
{4\gamma}{r^{2}}v(r)\right)  w(\theta)\nonumber\\
&  =\Lambda v(r)w(\theta)\ \text{for all}\ r\in(0,\infty)\ \text{and
all}\ \theta\in\lbrack0,\pi]. \label{eq16.7}%
\end{align}

For (\ref{eq16.7}) to hold requires that the function $v(\cdot),$ on gathering
up terms, has to satisfy the ordinary differential equation, see \cite[Section
1, (4)]{L1}, for all $r\in(0,\infty),$%
\begin{equation}
v^{(4)}(r)+\dfrac{2}{r}v^{(3)}(r)+\left(  -\dfrac{9}{r^{2}}-\gamma\right)
v^{\prime\prime}(r)+\left(  \dfrac{9}{r^{3}}-\dfrac{\gamma}{r}\right)
v^{\prime}(r)-\Lambda v(r)=0. \label{eq16.8}%
\end{equation}
This last equation may be written in the Lagrange symmetric form%
\begin{equation}
(rv^{\prime\prime}(r))^{\prime\prime}-((9r^{-1}+\gamma r)v^{\prime
}(r))^{\prime}=\Lambda rv(r)\ \text{for all}\ r\in(0,\infty), \label{eq16.9}%
\end{equation}
which is the Bessel fourth-order differential equation (\ref{eq3.1}) when
$\gamma=8M^{-1}.$

Thus separated solutions of the partial differential equation (\ref{eq16.2})
can be written in the form%
\begin{equation}
u(r,\theta)=v(r)w(\theta)\ \text{for all}\ r\in(0,\infty)\ \text{and}%
\ \theta\in\lbrack0,2\pi], \label{eq16.10}%
\end{equation}
where $w(\cdot)$ is any trigonometrical solution (\ref{eq16.6}) of
(\ref{eq16.5}), and $v(\cdot)$ is any solution of the fourth-order Bessel
equation (\ref{eq16.9}) for any choice of the spectral parameter $\Lambda.$

Defining the partial differential expression $P_{\gamma}[\cdot],$ for
$\gamma\in(0,\infty),$ by%
\begin{equation}
P_{\gamma}[u]:=\left(  \nabla^{4}u-\gamma\nabla^{2}u-\frac{4\gamma}{r^{2}%
}u\right)  \label{eq16.11}%
\end{equation}
it is shown in \cite{L1} that $P_{\gamma}$ is a formally symmetric linear
partial differential expression in $L^{2}(\mathbb{E}^{2})$, using polar
co-ordinates $(r,\theta).$

Some early studies indicate that there may be problems in applied mathematics,
for which the partial differential equation $P_{\gamma}[u]=\Lambda u$ is
involved in one or more of the associated mathematical models.

\end{document}